\documentclass[10pt]{amsart}
\usepackage{amsmath,mathrsfs,cite}

\usepackage{graphicx,wasysym}
\usepackage[small]{caption}
\usepackage[colorlinks=true,urlcolor=blue,
citecolor=red,linkcolor=blue,linktocpage,pdfpagelabels,bookmarksnumbered,bookmarksopen]{hyperref}
\usepackage[english]{babel}
\usepackage[utf8]{inputenc}
\usepackage{enumerate}
\usepackage[hyperpageref]{backref}

\usepackage[left=2.8cm,right=2.8cm,top=2.7cm,bottom=2.7cm]{geometry}

\newcommand{\R}{{\mathbb R}}

\newcommand {\Div} {\mbox{\rm div}} 

\numberwithin{equation}{section}

\theoremstyle{definition}

\newcommand{\bgset}[1]{\big\{#1\big\}}

\newcommand{\norm}[2][]{\left\|#2\right\|_{#1}}

\DeclareMathOperator{\dvg}{div}

\title[First eigenpair for problems with variable exponents]{Computing 
the first eigenpair for \\ problems with variable exponents}

\author{Marcello Bellomi}
\author{Marco Caliari}
\author{Marco Squassina}

\address{Dipartimento di Informatica,
Universit\`a degli Studi di Verona
\newline\indent
C\'a Vignal 2, Strada Le Grazie 15,
I-37134 Verona, Italy}
\email{marco.caliari@univr.it}
\email{marco.squassina@univr.it}

\thanks{The research of the third author was partially
supported by the 2009 italian PRIN project: {\em Metodi Variazionali e Topologici
nello Studio di Fenomeni non Lineari}}

\begin{document}
	



\begin{abstract}
We compute the first eigenpair for variable exponent eigenvalue problems. We compare
the homogeneous definition of first eigenvalue with previous nonhomogeneous notions
in the literature. We highlight the symmetry breaking phenomena.
\end{abstract}

\maketitle

\section{Introduction}

The aim of this paper is the study of numerical solutions to the minimization problem introduced in \cite{Franzi} 
\begin{equation}
\label{Ray2}
 \lambda_1=\inf_{\overset{u\in W^{1,p(x)}_0(\Omega)}{u\neq 0}}\frac{\lVert \nabla u \rVert_{p(x)}}{\lVert u \rVert_{p(x)}}.
\end{equation}
Here $\Omega$ is a bounded domain in $\mathbb{R}^n$ and the variable 
exponent $p:\bar\Omega\to\R^+$ is a smooth function
such that $1<p^- \leq p(x) \leq p^+ < \infty$ for every $x\in\Omega$. 
The norm $\|\cdot\|_{p(x)}$ is the so-called Luxemburg norm
\begin{equation} \label{Luxdef}
\lVert f \rVert_{p(x)} = 
\inf \Big\{ \gamma > 0 : \int_{\Omega} \Big| \frac{f(x)}{\gamma} \Big|^{p(x)} \frac{1}{p(x)} \leq 1 \Big\}.
\end{equation}
If $p$ is a constant function, the problem reduces (up to a power $p$)
to the minimization of the quotient
\begin{equation}
\label{RayConstant}
 \inf_{\overset{u\in W^{1,p}_0(\Omega)}{u \neq 0}} \frac{\displaystyle\int_{\Omega} \lvert \nabla u \rvert^{p} }
 {\displaystyle\int_{\Omega} \lvert u \rvert^{p} }
\end{equation}
and, as it is known, its associated Euler--Lagrange equation is 
\begin{equation}
\label{ELConstant}
-\Div ( | \nabla u |^{p-2} \nabla u )=\lambda \left| u \right|^{p-2}u. 
\end{equation}
In this case,  we refer the reader to \cite{Lindqvist,stability} 
for the theoretical aspects and to \cite{computazionale} for a recent numerical
analysis.
The special case $p=2$ is the classical eigenvalue problem for the Laplacian 
$ -\Delta u=\lambda u,$
for which we refer the reader to \cite{tuttosul2}.
In general, in these type of problems, is crucial that some homogeneity holds, namely, 
if $u$ is a minimizer, so is $\omega u$ for any non-zero real constant $\omega$.
On the contrary, the quotient
\begin{equation}
\label{RayVariable}
 \frac{\displaystyle\int_{\Omega} \lvert \nabla u \rvert^{p(x)} }{\displaystyle\int_{\Omega} \lvert u \rvert^{p(x)} }
\end{equation}
with variable exponents fails to possess this feature. 
Therefore, as we point out in the following section,
its infimum over nontrivial functions of $W^{1,p(x)}_0(\Omega)$ turns out to be often equal 
zero and no minimizer exists \cite{cinese,sticinesi}. 
A way to avoid this collapse is to impose the constraint $ \int_{\Omega} \lvert u \rvert^{p(x)} \textrm{d}x =C.$
Unfortunately, doing so, minimizers obtained for different normalization 
constants $C$ are difficult to compare. 
For a suitable $p(x)$, it could even happen that any $\lambda>0$ is an eigenvalue for some choice of $C$. Thus \eqref{RayVariable} is not a proper generalization of \eqref{RayConstant}, which has  
well defined (variational) eigenvalues, although the full spectrum is not completely understood yet.
A way to avoid this situation is to use the Rayleigh quotient \eqref{Ray2}, restoring 
the necessary homogeneity. In the integrand of \eqref{Luxdef}, the use of the measure $p(x)^{-1} \textrm{d}x$ simplifies the Euler-Lagrange equation.
The Sobolev inequality \cite{Base} $\| v \|_{p(x)} \leq 
C\|\nabla v \|_{p(x)},$  
shows that $\lambda_1 > 0$. It is easy to see that \eqref{Ray2} has a non-negative minimizer. 
Pick a minimizing sequence of $v_j$, namely $\lVert v_j \rVert_{p(x)} =1$ and 
$\lVert \nabla v_j \rVert_{p(x)}=\lambda_1+o(1).$ By Rellich theorem for variable 
Sobolev exponents \cite{Base}, up to a subsequence, we find $u$ such that 
$v_{j} \rightarrow u$ in $L^{p(x)} (\Omega)$ and $\nabla v_{j} \rightharpoonup \nabla u$ in $L^{p(x)} (\Omega)$. This yields
$\lambda_1\leq \lVert \nabla u \rVert_{p(x)}/\lVert u \rVert_{p(x)} \leq \lim_{j \to \infty} \lVert \nabla v_{j} \rVert_{p(x)}/
\lVert v_{j} \rVert_{p(x)} = \lambda_1.$
Notice that if $u$ is a minimizer so is $|u|\geq 0$.
By the maximum principle of \cite{harjulehto}, $u$ has a fixed sign.
In \cite{Franzi} the Euler-Lagrange equation for a minimizer is derived.
Precisely, it holds
\begin{equation}
\label{euler-lag-variab}
\int_{\Omega} \left| \frac{\nabla u}{K} \right|^{p(x)-2} \langle  \frac{\nabla u}{K} , \nabla \eta \rangle = \lambda_1 S \int_{\Omega} \left| \frac{ u}{k} \right|^{p(x)-2} \frac{u}{k} \eta,
\,\,\quad \forall\eta \in C^{\infty}_0 (\Omega),
\end{equation} 
where we have set 
$$
K= \lVert \nabla u \rVert_{p(x)},\quad\,\, 
k=\lVert u \rVert_{p(x)},\quad\,\,
\lambda_1 = \frac{K}{k},\quad\,\,
S=
{\Big(\displaystyle\int_{\Omega} \left| \frac{u}{k} \right|^{p(x)} \Big)^{-1}}
\displaystyle\int_{\Omega} \left| \frac{ \nabla u}{K} \right|^{p(x)}\!.
$$
More generally, $\lambda\in\R$ is eigenvalue if there exists $u \in W^{1,p(x)}_0 (\Omega)$, $u \not \equiv 0$, such that 
\begin{equation} \label{basic}
\int_{\Omega} \left| \frac{\nabla u}{K} \right|^{p(x)-2} \langle  \frac{\nabla u}{K} , \nabla \eta \rangle 
= \lambda S \int_{\Omega} \left| \frac{ u}{k} \right|^{p(x)-2} \frac{u}{k} \eta, 
\quad\,\, \forall\eta \in C^{\infty}_0 (\Omega).
\end{equation}
It follows from the regularity theory developed in \cite{acerbming} that the solutions to \eqref{basic}
are continuous provided that $p(x)$ is H\"older continuous.
If $\lambda_1$ is the minimum in \eqref{Ray2}, we have 
$\lambda \geq \lambda_1$ in \eqref{basic}, thus $\lambda_1$ is the first eigenvalue and a corresponding solution is the first eigenfunction. 
Contrary to the constant exponent case \cite{Lindqvist,stability}, it is currently 
unknown if, in the variable exponent case, the first eigenvalue is simple, and 
if a given positive eigenfunction is automatically a first one.
Concerning higher eigenvalues, in \cite{persqu} the authors have recently proved that
there is a sequence of eigenvalues of \eqref{basic} with $\lambda_j \nearrow \infty$ and if 
\[
\sigma = n \left(\frac{1}{p^-} - \frac{1}{p^+}\right) < 1,
\]
then there are constants $C_1, C_2 > 0$, that depend only on $n$ and $p^\pm$, such that
\begin{equation}
\label{eigenvestim}
C_1\, |\Omega|\, \lambda^{n/(1 + \sigma)} \le \# \bgset{j : \lambda_j < \lambda} \le C_2\, |\Omega|\, \lambda^{n/(1 - \sigma)} \qquad \text{for $\lambda > 0$ large},
\end{equation}
where $|\Omega|$ is the Lebesgue measure of $\Omega$. Observe that, in the case of constant $p$, \eqref{basic} reduces not exactly
to \eqref{ELConstant}, which is homogeneous of degree $p-1$, but rather to the problem (homogeneous of degree $0$)
\[
- \dvg \Big(\frac{|\nabla u|^{p-2}\, \nabla u}{\norm[p]{\nabla u}^{p-1}}\Big) = \lambda\, \frac{|u|^{p-2}\, u}{\norm[p]{u}^{p-1}}, \quad\,\, u \in W^{1,p}_0(\Omega)
\]
Thus \eqref{eigenvestim}
should be compared to $C_1\, |\Omega|\, \lambda^{n/p} \le \# \bgset{j : \lambda_j < \lambda} \le C_2\, |\Omega|\, \lambda^{n/p}$, obtained in \cite{MR1017063}. 

\subsection{A different notion in the literature}

We compare the minimization procedure
with the Rayleigh quotient with Luxemburg norm 
and that without it, namely
$$
\inf_{\overset{u\in W^{1,p(x)}_0(\Omega)}{u\neq 0}}
\frac{\displaystyle\int_{\Omega} |\nabla u|^{p(x)}}{
\displaystyle\int_{\Omega} |u|^{p(x)}}. 
$$
In this framework, if $\lambda \in \mathbb{R}$ and 
$u \in W^{1,p(x)}_0 (\Omega)$ then $(u,\lambda)$ is called eigenpair if $u\neq 0$ and
$$
\int_{\Omega} |\nabla u|^{p(x)-2} \langle\nabla u, \nabla \eta \rangle= \lambda \int_{\Omega} | u|^{p(x)-2} u \eta , \quad\,\, 
\forall \eta \in W^{1,p(x)}_0 (\Omega). 
$$
Set $\Lambda =  \lbrace \lambda>0 : \lambda \textrm{ is an eigenvalue}  \rbrace.$
It is well known \cite{Lindqvist,stability}  that, if the function $p(x)$ is constant, then the problem has a sequence of eigenvalues, $\sup \Lambda = + \infty$ and $\inf \Lambda>0$.
In the general case, it follows from \cite{sticinesi} that
$\Lambda$ is a nonempty infinite set and $\sup \Lambda = + \infty$.
Define $\lambda_* := \inf \Lambda.$
We recall that we often have $\lambda_* =0$ (recall that $\lambda_1 > 0$ in~\eqref{Ray2}).
Consider the following Rayleigh quotients 
\begin{align*}
\mu_*  = \inf_{\overset{u \in W^{1,p(x)}_0 (\Omega)}{u\neq 0}} \frac{\displaystyle\int_{\Omega} \frac{|\nabla u|^{p(x)}}{p(x)} }{\displaystyle
\int_{\Omega} \frac{|u|^{p(x)}}{p(x)}  }, \qquad 
\overline{\mu}_* = \inf_{\overset{u \in W^{1,p(x)}_0 (\Omega)}{u\neq 0}} 
\frac{\displaystyle\int_{\Omega} |\nabla u|^{p(x)} }{\displaystyle
\int_{\Omega} | u|^{p(x)} }.
\end{align*}
Then, in \cite{sticinesi}, the authors prove that
$\lambda_* > 0 \Leftrightarrow \mu_* > 0 \Leftrightarrow \overline{\mu}_* > 0$.
Furthermore, if there is an open subset $U \subset \Omega$ and a point $x_0 \in U$ 
such that $p(x_0) < p(x)$ (or $>$) for all $x \in \partial U$, then $\lambda_* =0$
\cite[Theorem 3.1]{sticinesi}.
In particular, if $p(x)$ has strictly local minimum (or maximum) points in $\Omega$, then $\lambda_* =0$. There are also statements giving some sufficient conditions for $\inf \Lambda > 0$. Let $n>1$. If there is a vector $\ell \in \mathbb{R}^n \setminus \lbrace 0 \rbrace$ such that, for any $x \in \Omega$, 
the map $t\mapsto 
p(x+t\ell)$ is monotone on $\{ t : x+t\ell \in \Omega\}$, then $\lambda_* >0$
\cite[Theorem 3.3]{sticinesi}.
If $n=1$, then $\lambda_* >0$ if and only if the function $p(x)$ is monotone 
\cite[Theorem 3.2]{sticinesi}.

\begin{figure}[!ht]
\centering
\hfill
\includegraphics[scale=0.47]{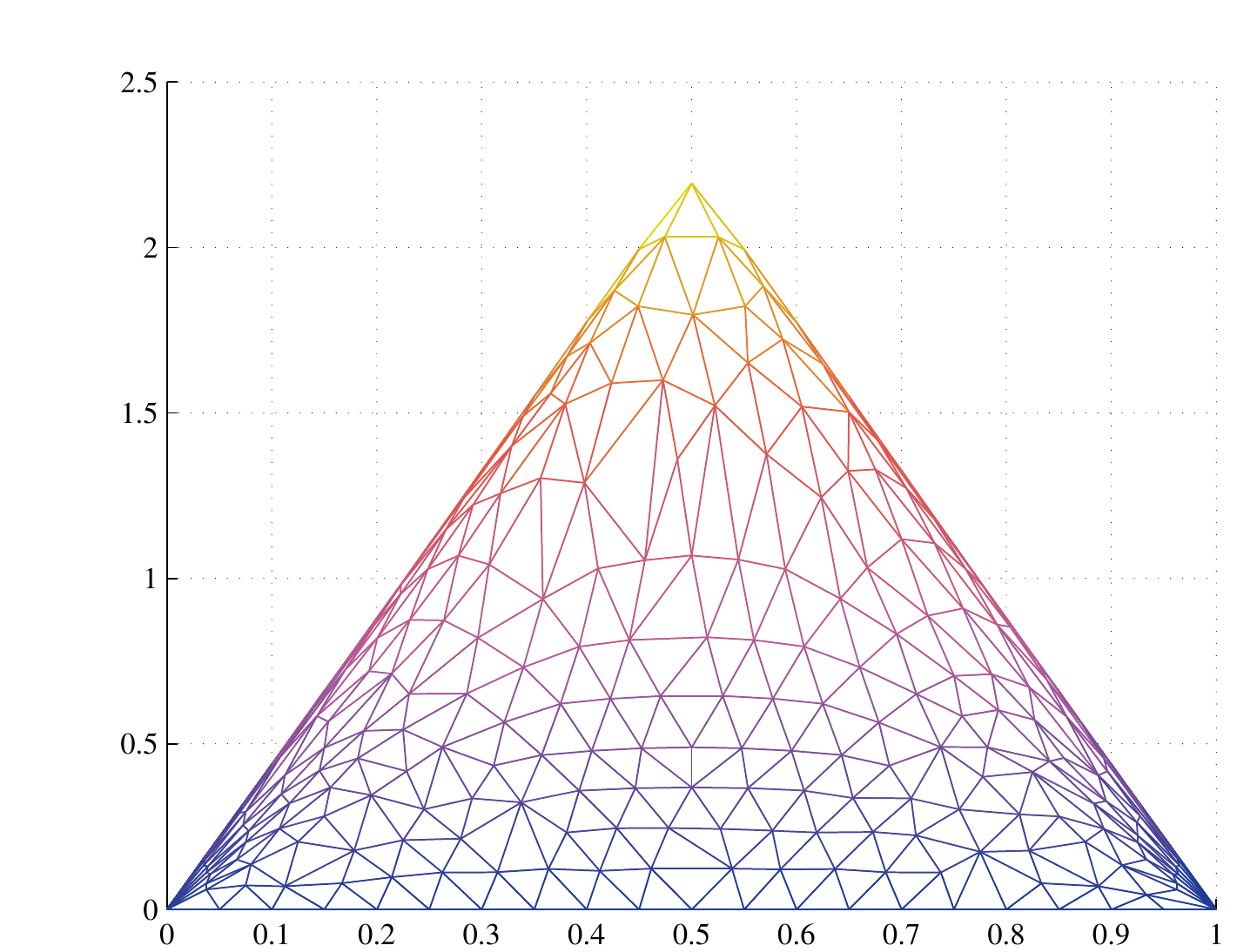}\hfill
\includegraphics[scale=0.47]{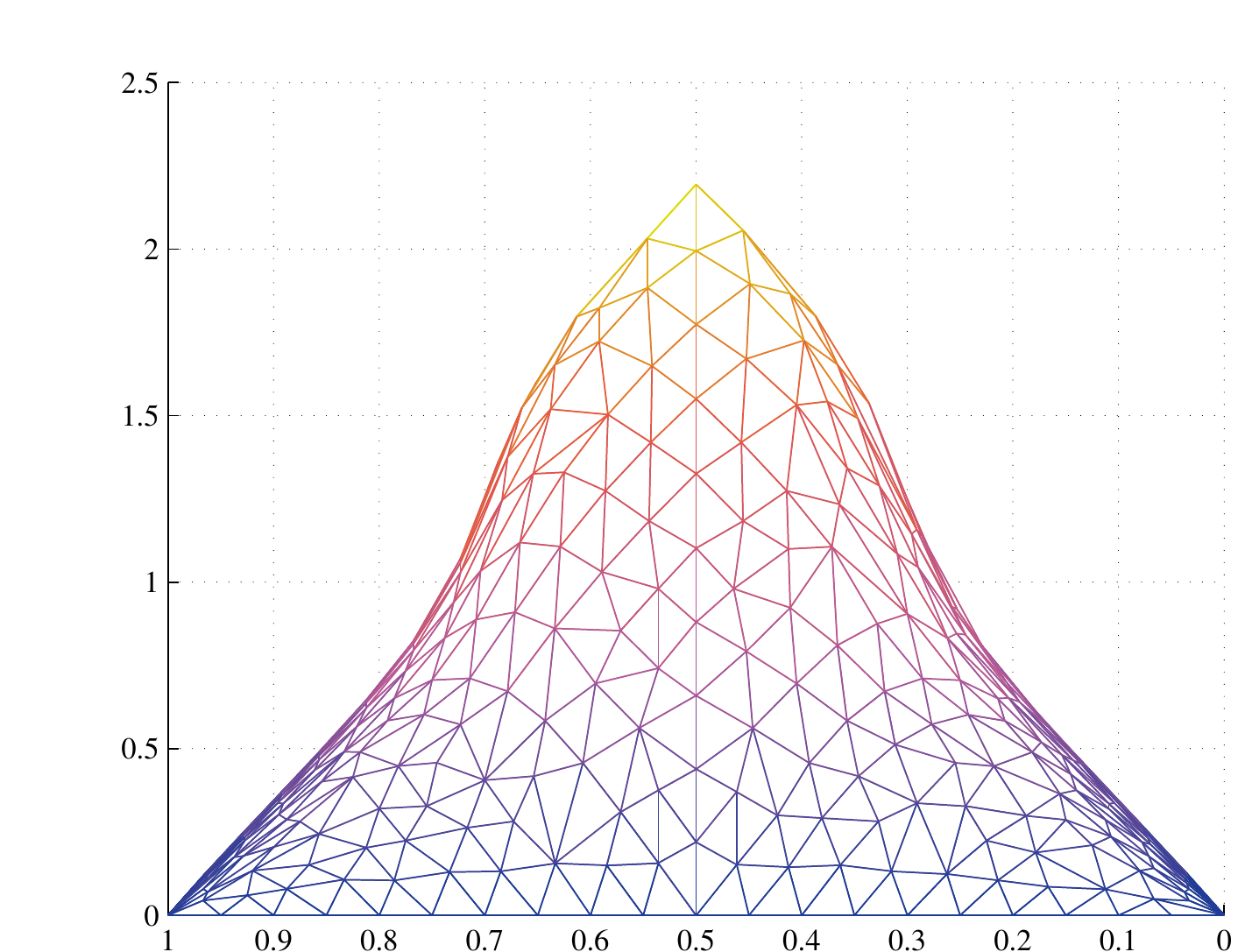}\hfill{}

\includegraphics[scale=0.47]{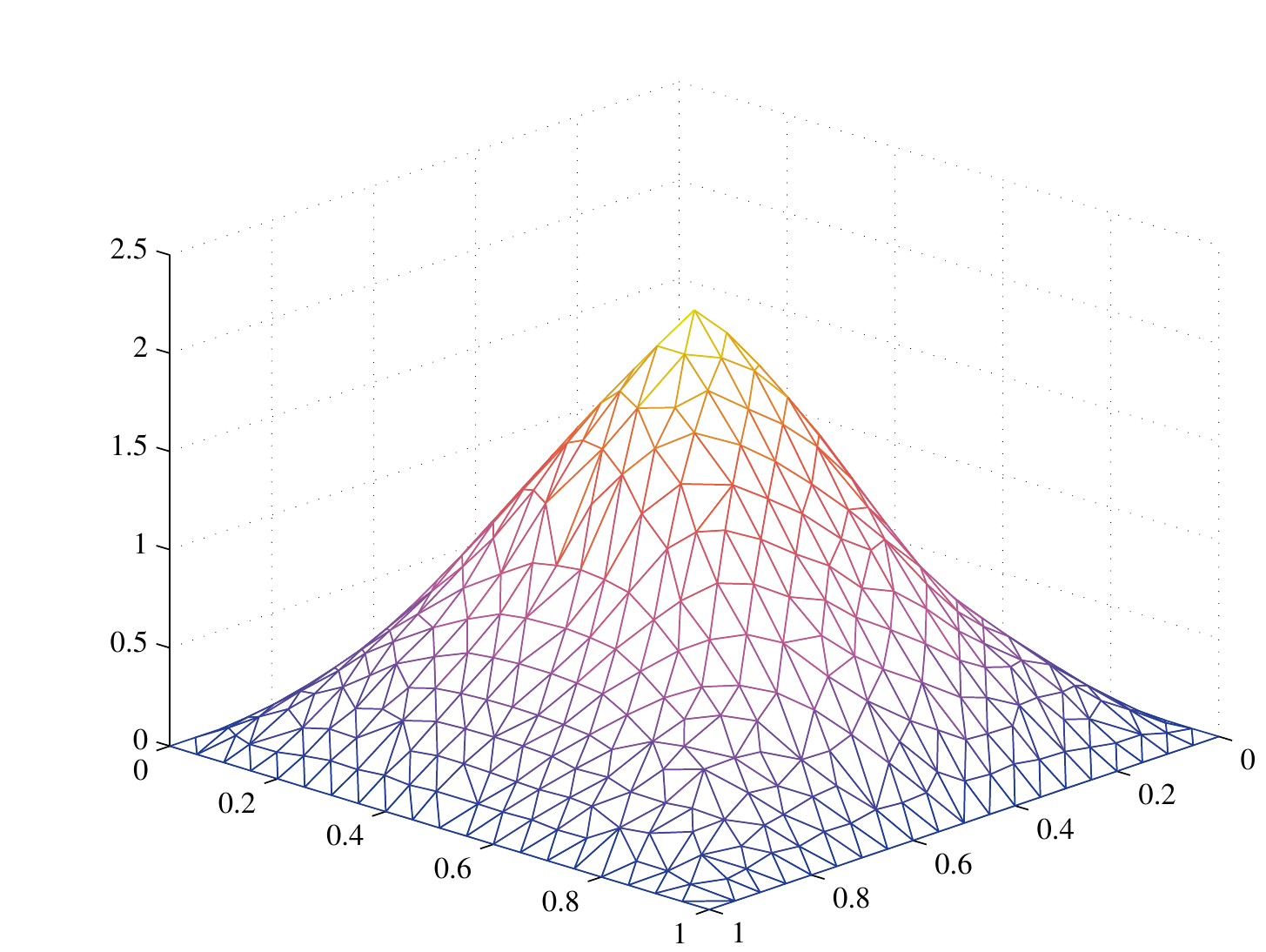}
\caption{Unit square, $p(x,y)=5+3\sin(3\pi x)$, $x$-axis (left), 
$y$-axis (right) and top views}
\label{fig:square}
\end{figure}

\begin{figure}[!ht]
\centering
\hfill
\includegraphics[scale=0.47]{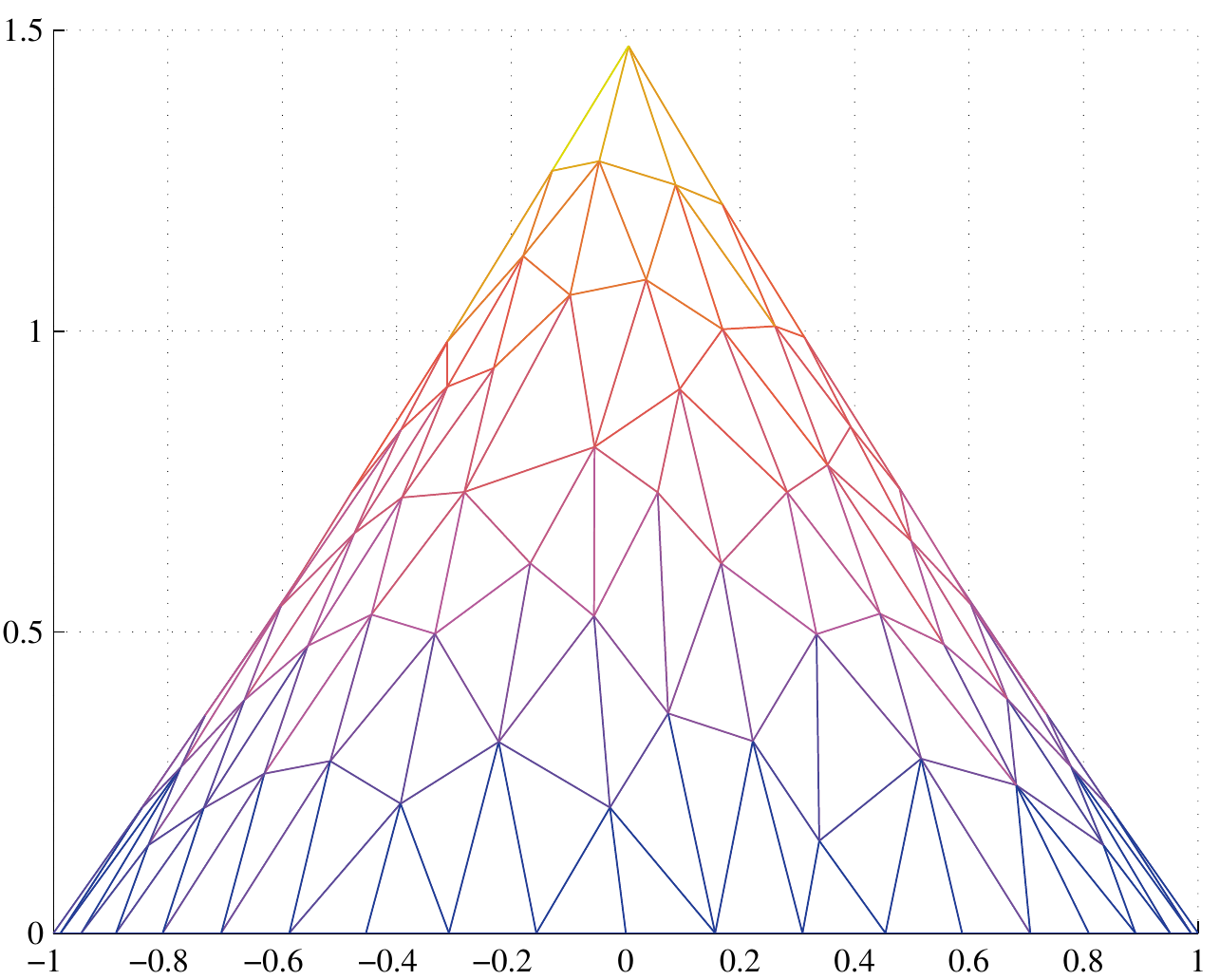}\hfill
\includegraphics[scale=0.47]{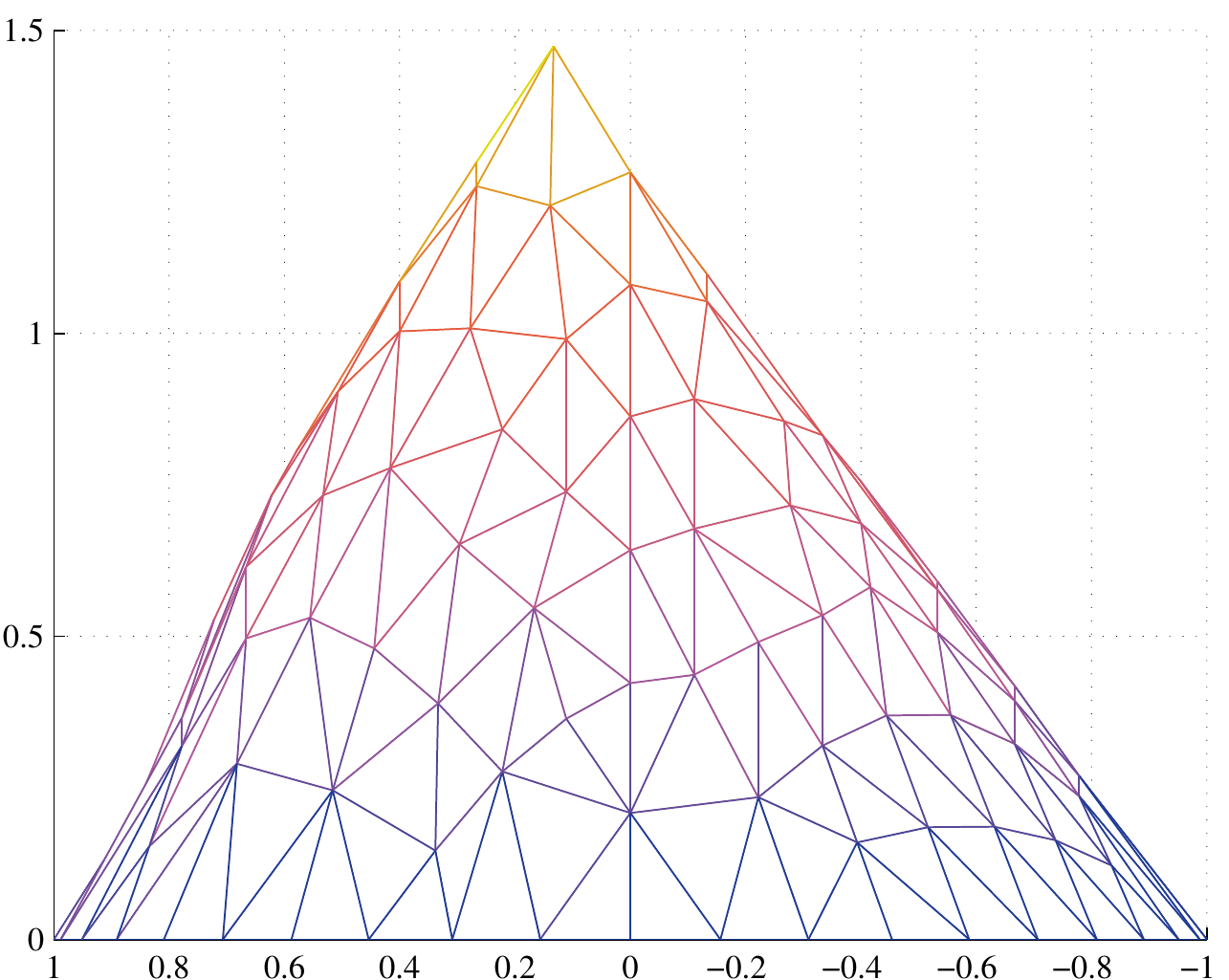}\hfill{}

\includegraphics[scale=0.47]{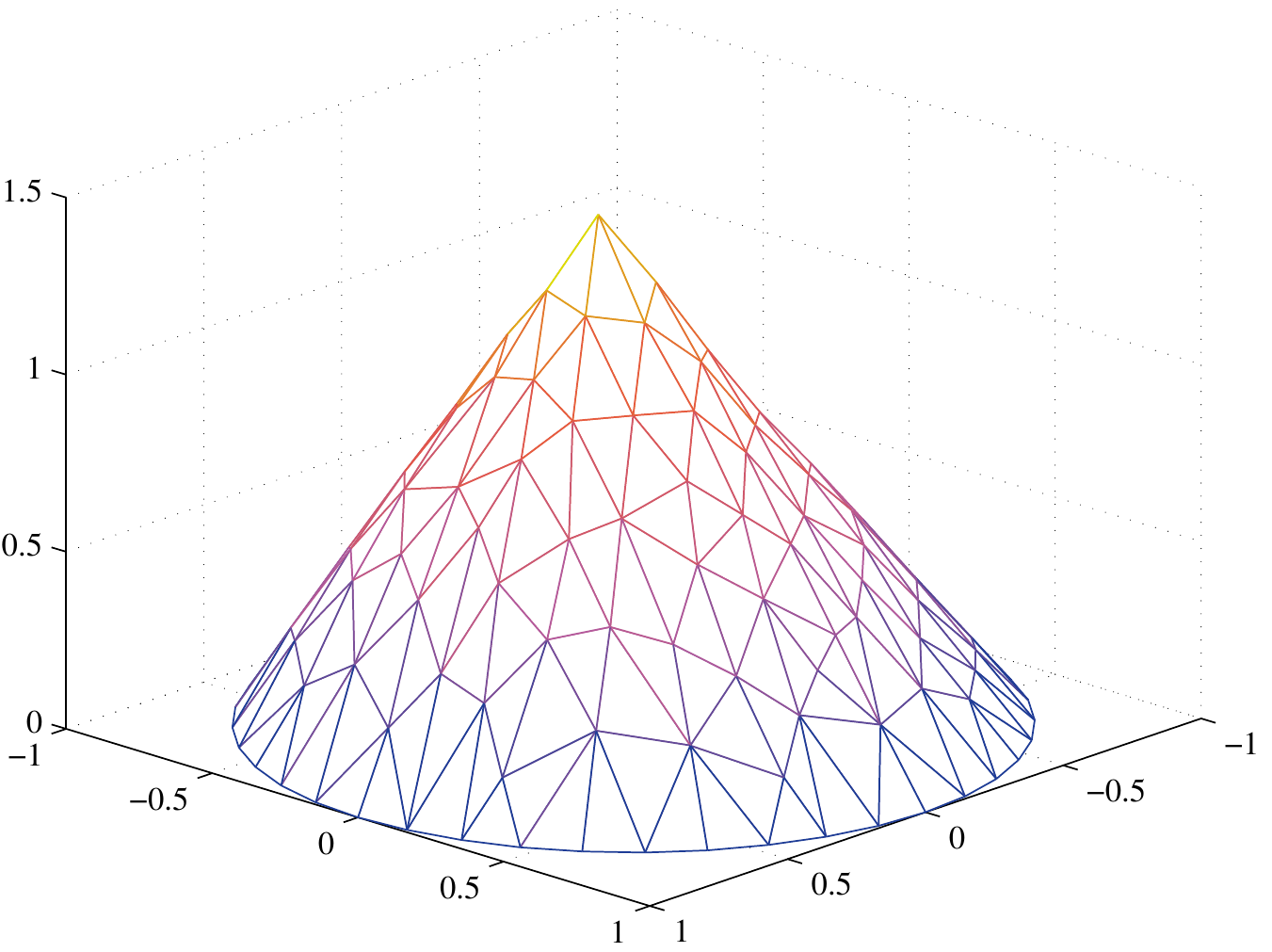}
\caption{Unit disk, $p(x,y) = 11+9\sin(2 \pi x)$, $x$-axis (left), 
$y$-axis (right) and top views}
\label{fig:disk}
\end{figure}

\begin{figure}[!ht]
\centering
\hfill
\includegraphics[scale=0.47]{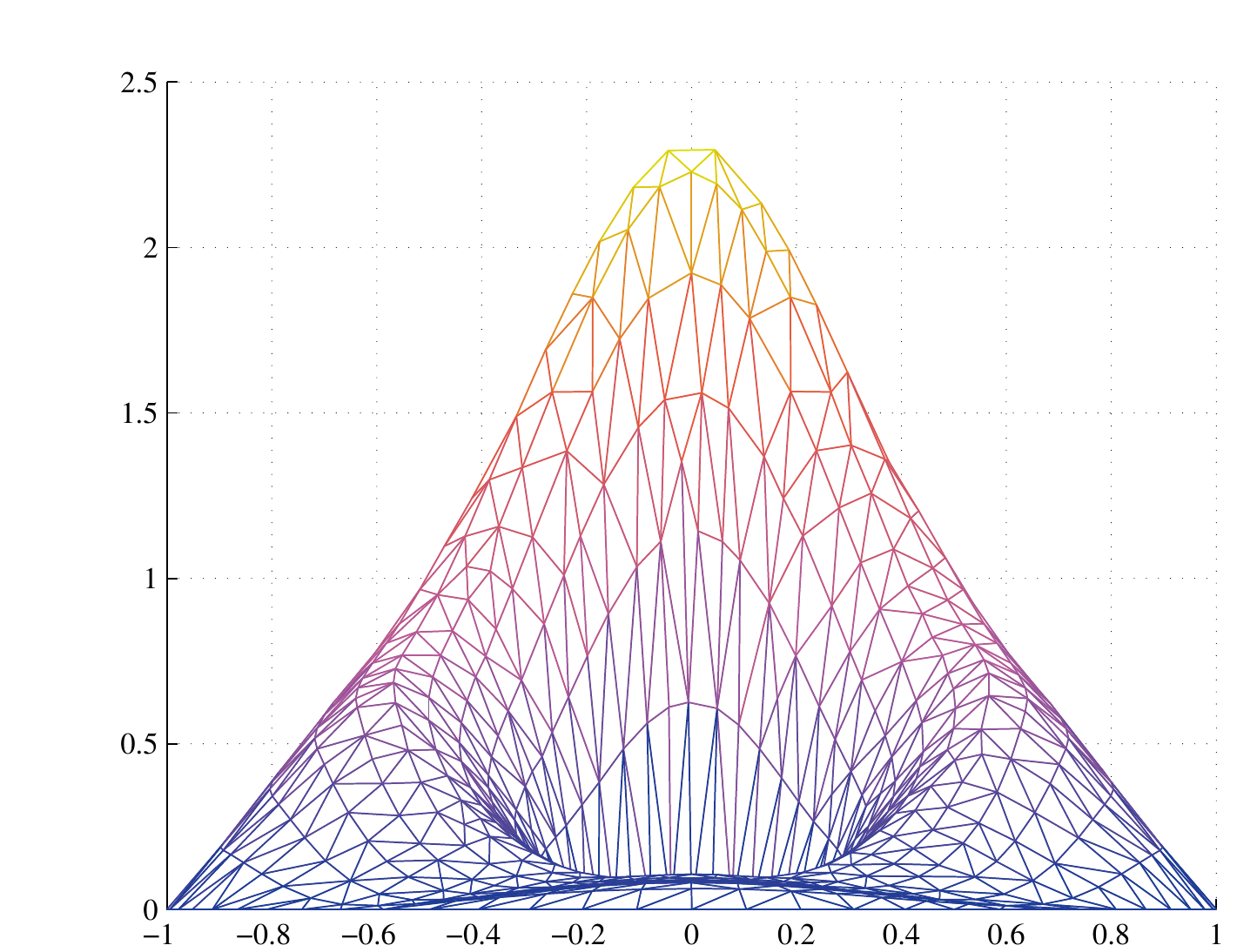}\hfill
\includegraphics[scale=0.47]{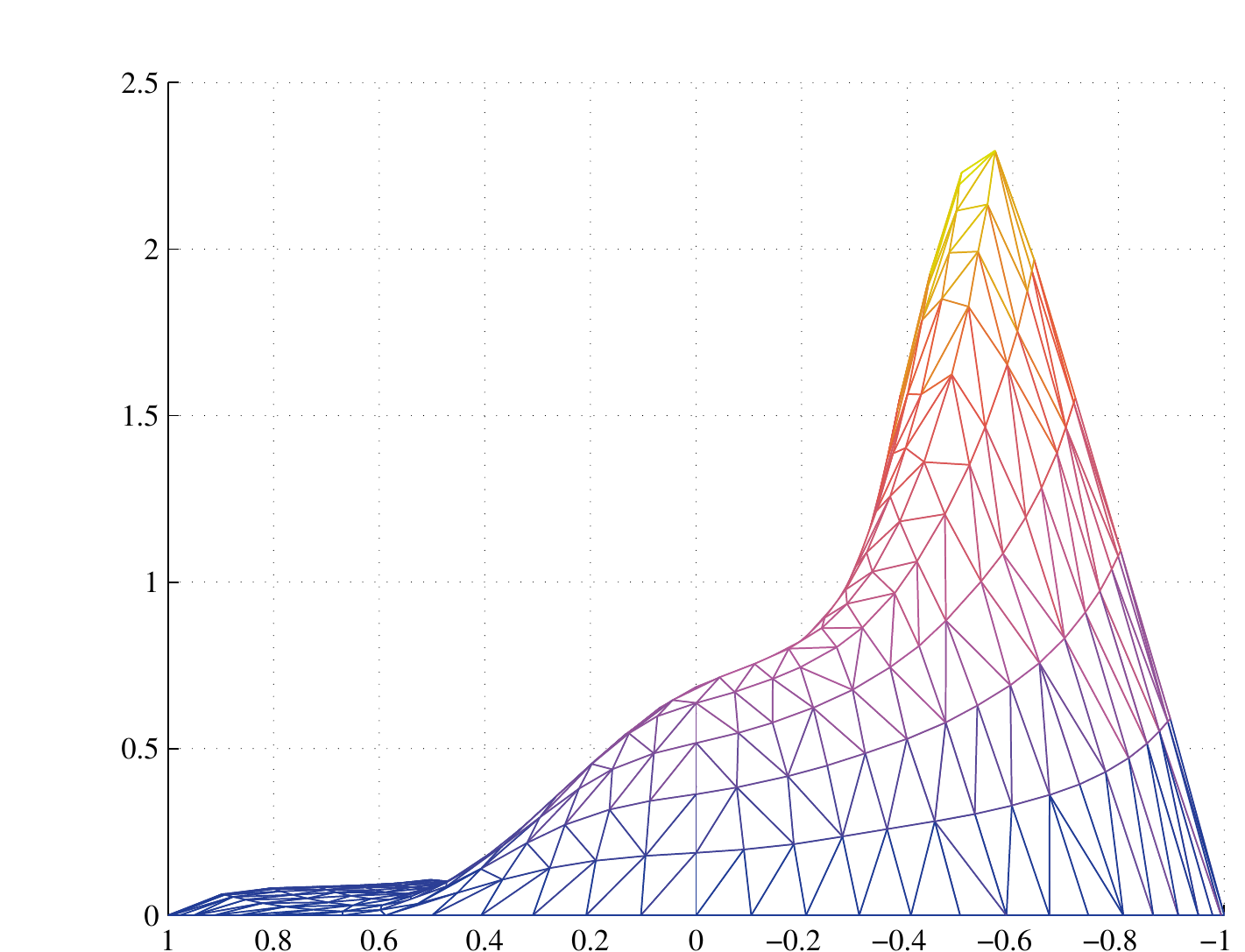}\hfill{}

\includegraphics[scale=0.47]{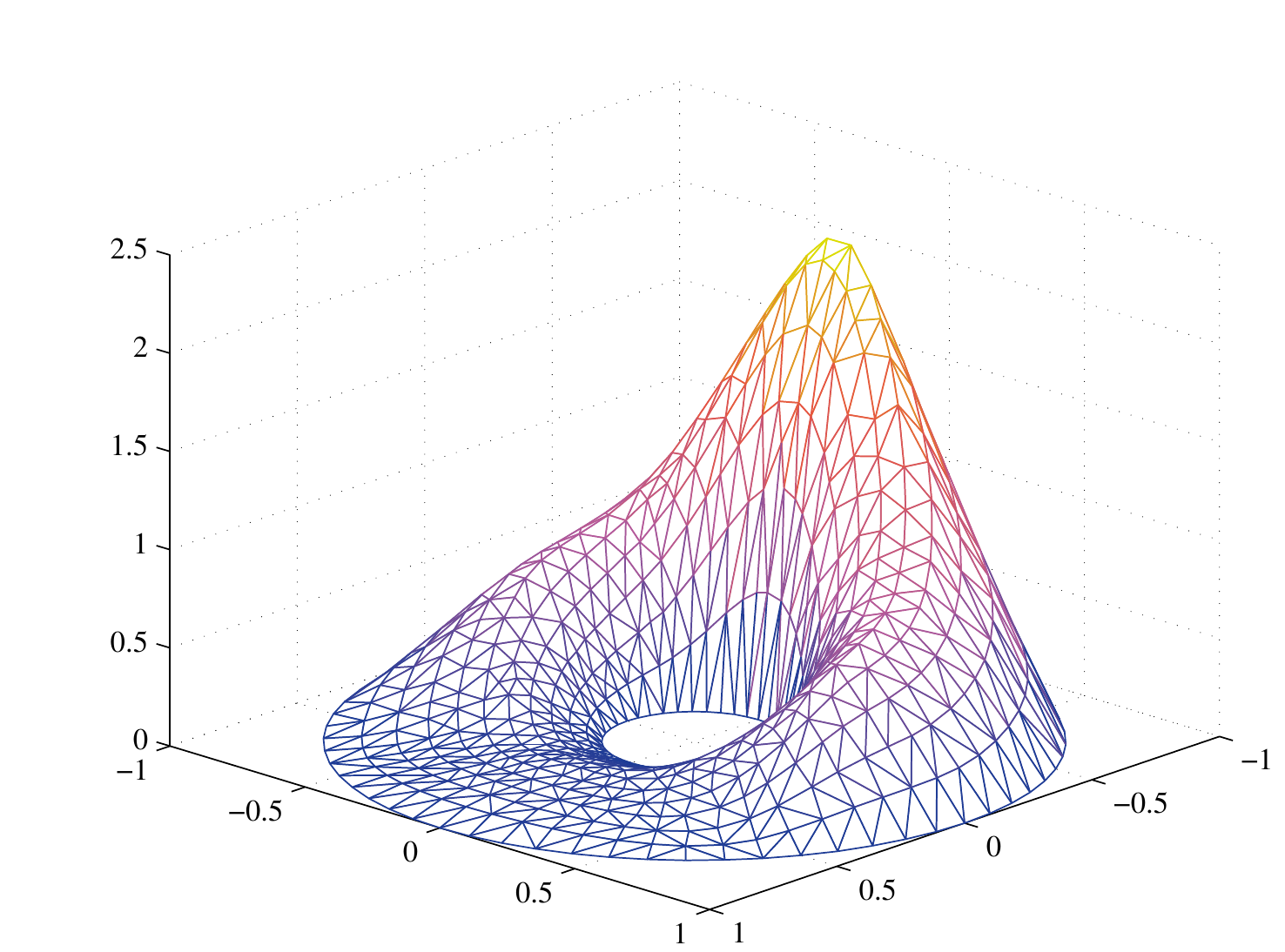}
\caption{Annulus, $p(x,y)=4+2\sin(2\pi x)$, $x$-axis (left), 
$y$-axis (right) and top views}
\label{fig:annulus}
\end{figure}

\section{An algorithm to compute the first eigenpair} \label{algo}
In this section we briefly describe an algorithm to approximate $\lambda_1$
in~\eqref{Ray2} and compute the corresponding eigenfunction. We start defining
\begin{equation*}
\Lambda_1:= \inf_{u \in W^{1,{p(x)}}_0(\Omega)\setminus \{0\}} \frac{\|\nabla u\|_{p(x)}^2}{\|u\|_{p(x)}^2}=
\inf_{u \in W^{1,{p(x)}}_0(\Omega)\setminus \{0\}} \frac{R(u)}{S(u)}=\lambda_1^2
\end{equation*}
It is now possible to apply the inverse power method, where the 
$(j+1)$-th iteration is
\begin{subequations}
\begin{align}
\tilde u^{j+1}&=\arg\min_u(R(u)-\nabla S(u^j)u)=
\arg\min_u J(u)\label{eq:inner}\\
u^{j+1}&=\frac{\tilde u^{j+1}}{S(\tilde u^{j+1})^{1/2}}\label{eq:normal}\\
\Lambda_1^{j+1}&=\frac{R(u^{j+1})}{S(u^{j+1})}
\end{align}
\end{subequations}
where $u^j$ is the result of the previous iteration and, by~\eqref{eq:normal},
has Luxemburg norm equal to 1.
It is possible to show (see~\cite{HB10}) that the algorithm converges
to a critical point of $R(u)/S(u)$, even if it is not possible in general to
prove convergence to the smallest eingevalue. However, a good choice of the 
initial
guess $u^0$ can reasonably assure that the result is the smallest eigenvalue.
For the computation of $\nabla S(u)\eta$, for given $u$ and $\eta$, in the
so called \emph{inner problem}~\eqref{eq:inner},
we observe that if $u\ne 0$ its
Luxemburg norm $\gamma(u)=\|u\|_{p(x)}$ is implicitly defined by
\begin{equation}\label{eq:implicit}
F(u,\gamma)=\int_\Omega\left|\frac{u(x)}{\gamma}\right|^{p(x)}
\frac{1}{p(x)}-1=0
\end{equation}
Therefore, we can use the differentiation of implicit functions to get
\begin{equation*}
\nabla \|u\|_{p(x)}\eta=-\frac{\nabla_u F(u,\gamma)\eta}{\partial_\gamma F(u,\gamma)}
=\frac{\displaystyle\int_\Omega\left|\frac{u}{\|u\|_{p(x)}}\right|^{p(x)-2}
\frac{u}{\|u\|_{p(x)}}\eta}{\displaystyle\int_\Omega\left|\frac{u}{\|u\|_{p(x)}}\right|^{p(x)}}
\end{equation*}
from which
\begin{equation*}
\nabla S(u)\eta=2\frac{\displaystyle\int_\Omega\left|\frac{u}{\|u\|_{p(x)}}\right|^{p(x)-2}
u\eta}{\displaystyle
\int_\Omega\left|\frac{u}{\|u\|_{p(x)}}\right|^{p(x)}}
\end{equation*}
Since we are mainly interested is some particular two-dimensional domains, 
such as a rectangle, a disk or an annulus, we approximated the problem
by the finite element method
which well adapts to different geometries by constructing an appropriate
discretization mesh. 
In particular, we used the tool FreeFEM++~\cite{freefem} which
can handle minimization problems as~\eqref{eq:inner} through the
function
\verb+NLCS+ (nonlinear conjugate gradient method, Fletcher--Reeves 
implementation). Such a function requires the application of the gradient of 
$J(u)$ to a test function $\eta$
\begin{equation*}
\nabla J(u)\eta=
\frac{\displaystyle\int_\Omega\left|\frac{\nabla u}{\|\nabla u\|_{p(x)}}\right|^{p(x)-2}\langle\nabla u,\nabla \eta\rangle}{\displaystyle\int_\Omega\left|\frac{\nabla u}{\|\nabla u\|_{p(x)}}\right|^{p(x)}}-
\frac{\displaystyle\int_\Omega\left|\frac{u^j}{\|u^j\|_{p(x)}}\right|^{p(x)-2}
u^j\eta }{\displaystyle\int_\Omega\left|\frac{u^j}{\|u^j\|_{p(x)}}
\right|^{p(x)}}
\end{equation*}
and an initial guess which, for the $(j+1)$-th iteration of the inverse
power method, is $u_j/\Lambda_1^j$.
The stopping criterion for the inverse power method is based on the
difference of two successive approximations of $\Lambda_1$.

\subsection{Some details of the algorithm}
The algorithm described above requires recurrent computations of the Luxemburg
norm of a function. For a given $u\ne 0$, it is the zero of the
function $F(u,\cdot)$ defined in~\eqref{eq:implicit}. This is 
a $C^2(0,+\infty)$ convex and monotonically decreasing function in $\gamma$,
with $\lim_{\gamma\to 0^+}F(u,\gamma)=+\infty$ and 
$\lim_{\gamma\to+\infty}F(u,\gamma)=-1$. Therefore, it is possible to
apply the quadratically convergent Newton's method in order to find its 
unique zero, starting
with an initial guess $\gamma_0$ on its left hand side (i.e., such that 
$F(u,\gamma_0)>0$).
As pointed out above, the inverse power method cannot guarantee the convergence to
the smallest eigenvalue and relative eigenfunction. 
It is very reasonable to expect that if the initial
guess $u^0$ for the method is close enough to the smallest eigenfunction, then
the algorithm will converge to it. For $p\equiv2$ the problem essentially reduces
to the Helmholtz equation
\begin{equation*}
-\Delta u=\lambda_1u,
\end{equation*}
for which the eigenfunctions are well known for the domains we have
in mind. Therefore, starting with $p\equiv2$ and $u^0$ the
eigenfunction corresponding to the smallest eigenvalue for the Helmholtz
equation we moved to the desired $p(x)$ through a standard
continuation technique.

We tested our algorithm using both linear and quadratic finite elements and
checked the correct order of convergence (two and three for the rectangle and
one and two for the circular domains,  respectively). Moreover,
we checked that the results with constant $p$ were consistent with those 
reported in~\cite{computazionale}.
The results in the next section were obtained with quadratic finite elements.
We observed convergence of our algorithm also for some cases with $p(x)<2$, 
for which the Hessian of $J(u)$ degenerates and the nonlinear conjugate 
gradient is not guaranteed to converge. In this case some authors add
a regularization parameter to the functional $J(u)$ (see~\cite{computazionale}).
We moreover observed sometimes slow convergence of the nonlinear conjugate 
gradient.
In this case, a more sophisticated method, using the Hessian of $J(u)$ or an
approximation of it, could be employed. Another possibility would be
to use a preconditioner, either based on a low order approximation
of $J(u)$ (see again~\cite{computazionale} and reference therein) or on
a linearized version of $J(u)$. The implementation of a more robust and
fast algorithm, on which we are currently working, 
is beyond the scopes of this paper.

\begin{figure}[!ht]
\centering
\includegraphics[scale=0.26]{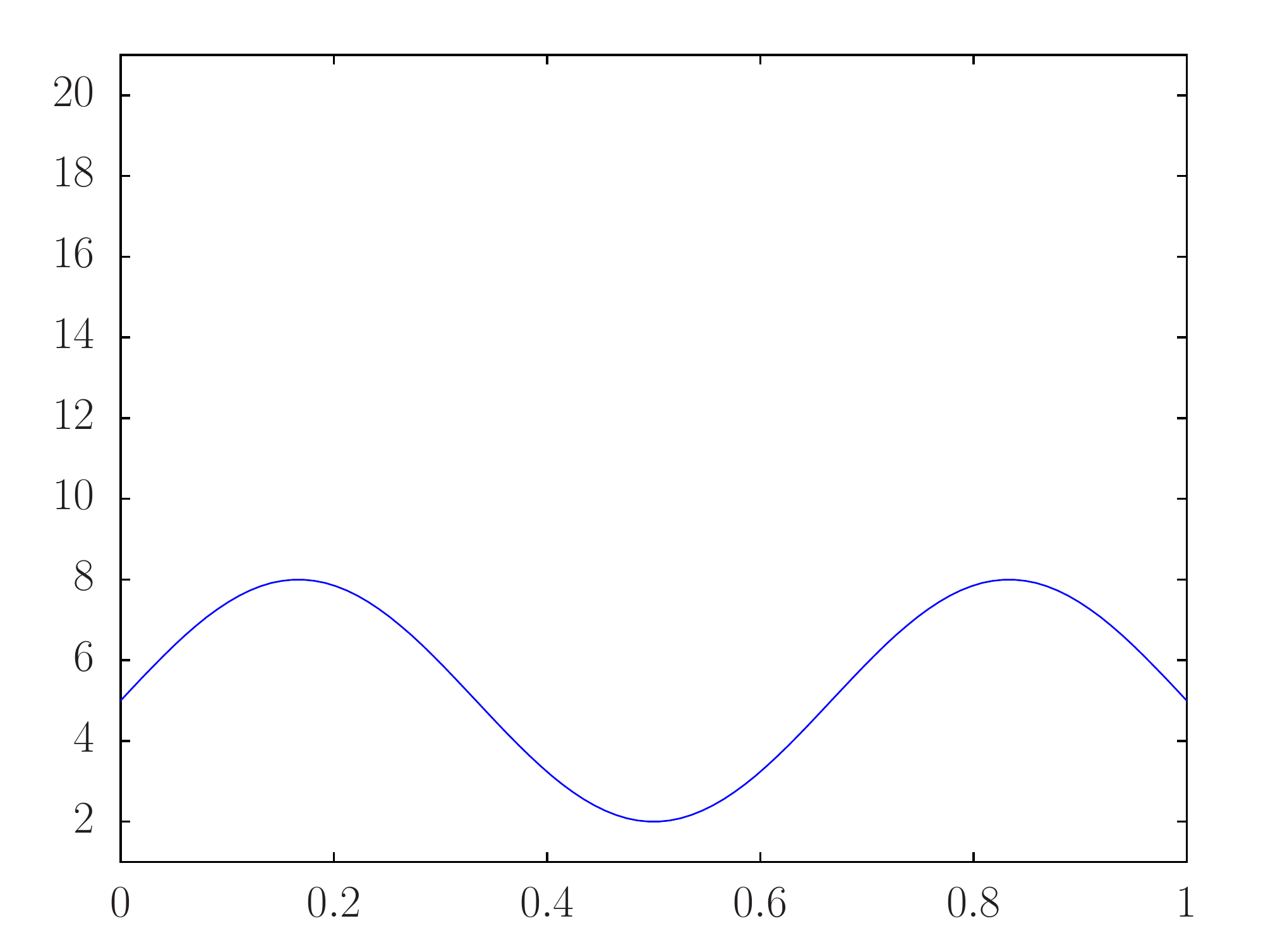}\hfill
\includegraphics[scale=0.26]{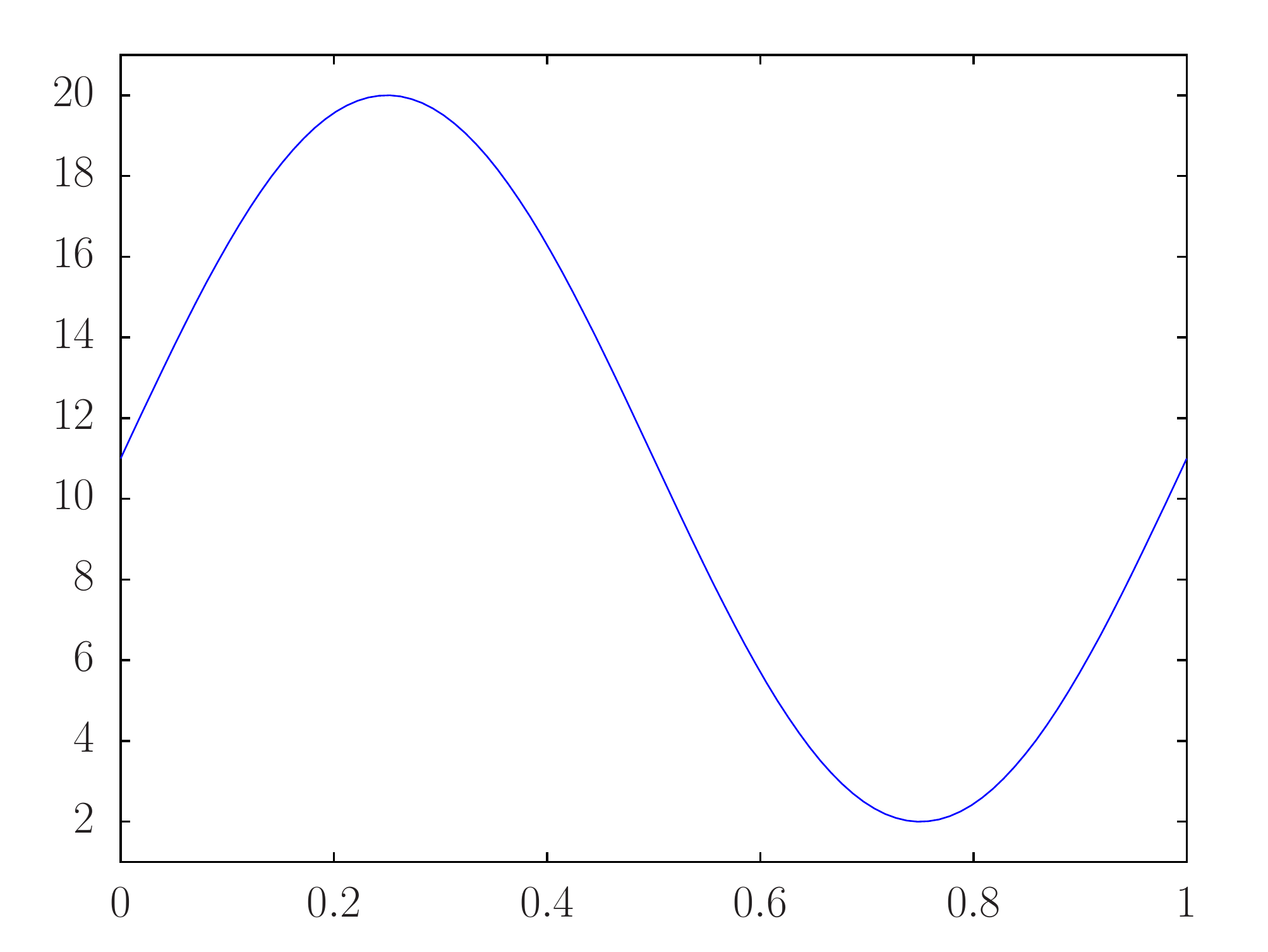}\hfill
\includegraphics[scale=0.26]{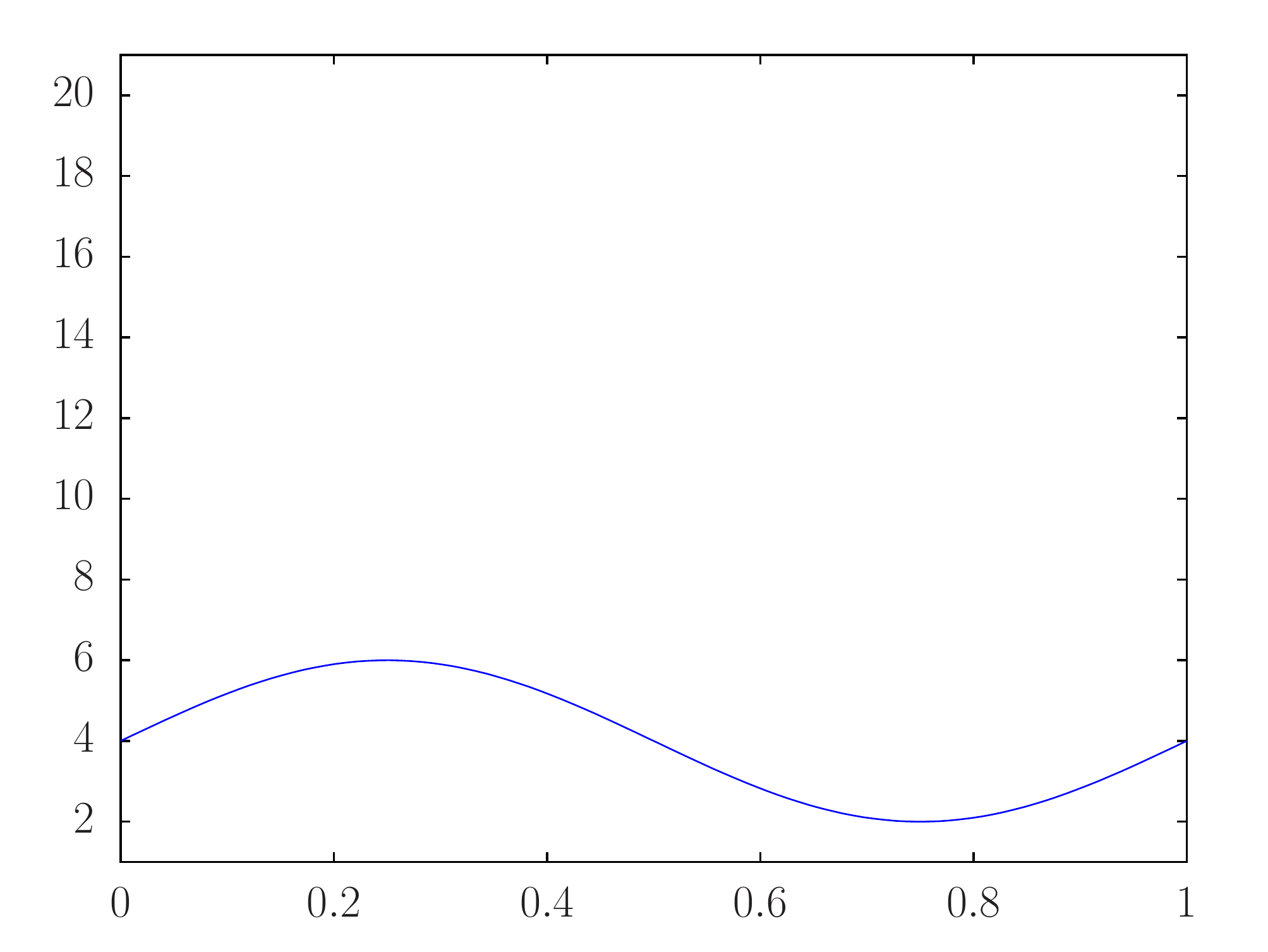}
\caption{$p(x,y)=5+3\sin(3\pi x)$,
$p(x,y) = 11+9\sin(2 \pi x)$ and
$p(x,y)=4+2\sin(2\pi x)$, respectively}
\label{fig:p-cases}
\end{figure}

\subsection{Examples and breaking symmetry}
We show in this section the results we obtained by applying the algorithm
described above to three test cases, in the square, the disk and the annulus,
respectively. For each test, we report the plot of the obtained eigenfunction
from the $x$-axis, the $y$-axis and from the top view, respectively.

\begin{figure}[!ht]
\centering
\hfill
\includegraphics[scale=0.35]{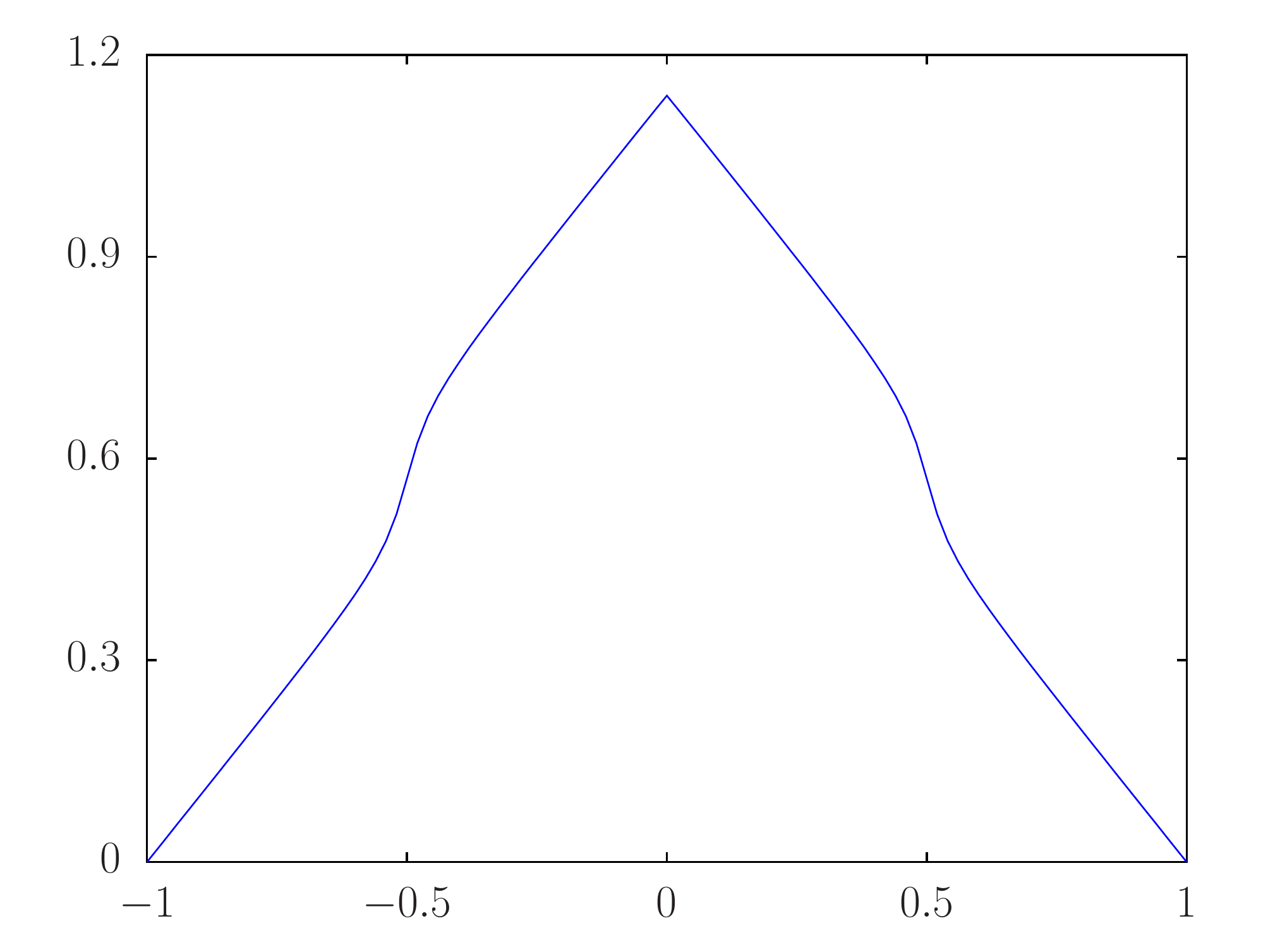}\hfill
\includegraphics[scale=0.35]{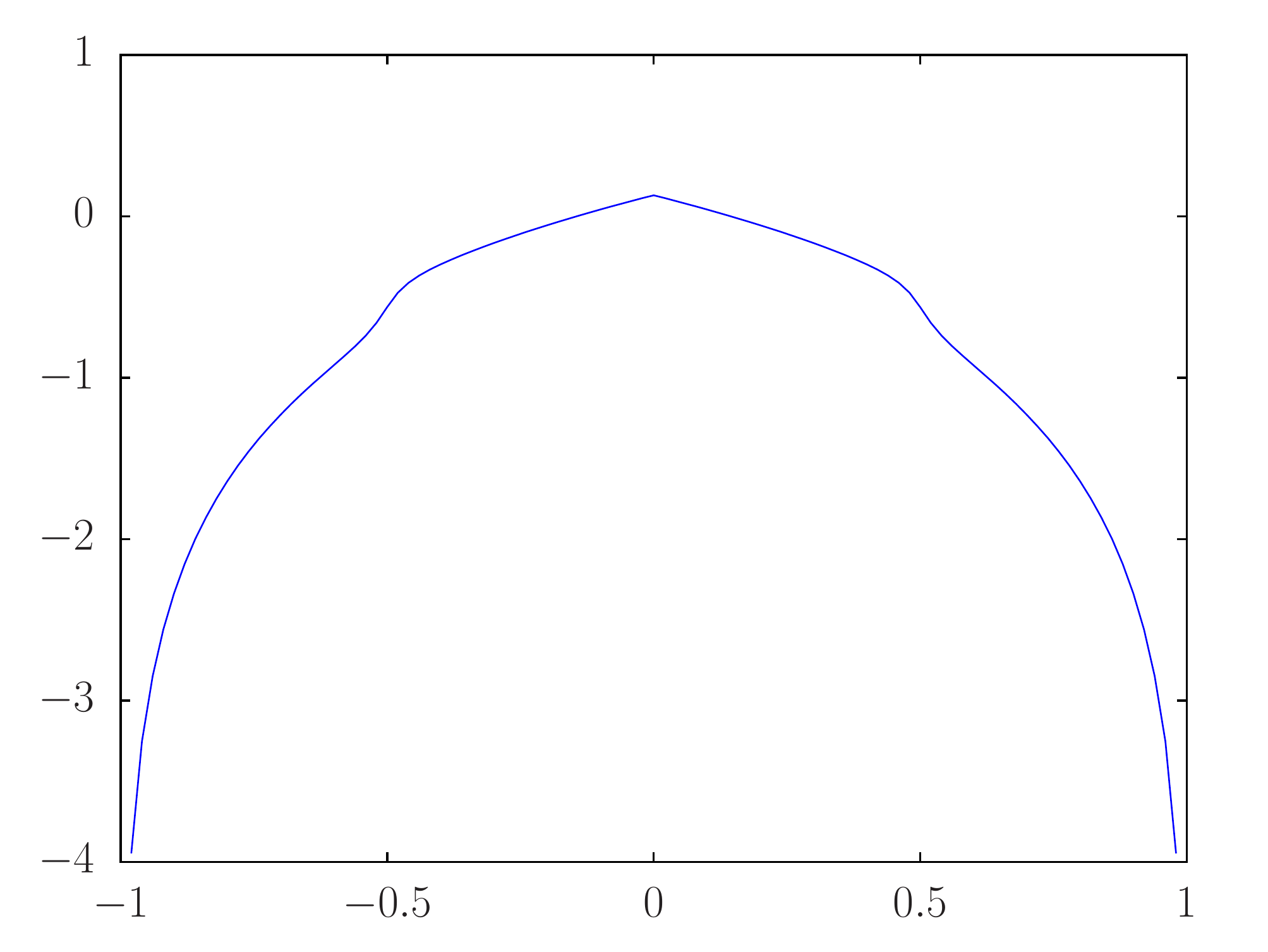}\hfill{}
\caption{Interval $[-1,1]$, $p(x)=28+26\cos(2\pi x)$, eigenfunction (left) and
its logarithm (right)}
\label{fig:1d}
\end{figure}

\noindent
The first case (Figure~\ref{fig:square}) refers to the {\em unit square} $[0,1]\times [0,1]$ 
and $p(x,y)=5+3\sin(3\pi x)$. The first plot (left) in this figure does not present special features because the exponent $p(x,y)$ 
is independent of $y$. In the second plot (right) the profile is different as it feels the diffusion variation in the $x$-variable. 
Both plots are symmetric with respect to the center of the domain, since
$p(x,y)$ has a center of symmetry in $(1/2, 1/2)$, see Figure \ref{fig:p-cases}.
\noindent
The second case (Figure~\ref{fig:disk}) refers to the {\em unit disk} with center
$0$ and radius $1$ and $p(x,y)=11+9\sin(2\pi x)$. The maximum of $p(x,y)$ is quite high (see Figure \ref{fig:p-cases}) 
and the profile is reminiscent of the one for the limiting
case $p= \infty$. The second plot is remarkable because we loose the symmetry for the center of the domain.
\noindent
The third case (Figure~\ref{fig:annulus}) refers to the {\em annulus} with center 0,
external radius $1$ and internal radius $0.25$ and $p(x,y)=4+2\sin(2\pi x)$. The resulting 
eigenfunction is more shifted than the case in the unit square and the case in the unit disc. 
Even the shape of the domain influences the contour, here we see 
that the annulus reflects intensely the mold of the exponent. 
The first plot (left) maintains the center of symmetry in $(0,0)$ and $p(x,y)$ does not depend on $y$.

\noindent
Already in the one dimensional case, it is evident that the the logarithm of 
first eigenfunction is not a concave function, in general, contrary to the constant exponent case
where this was proved to be true \cite{sakag}. See Figure \ref{fig:1d} for an example in this regard.

\noindent
In the case where $\Omega$ has some symmetry and $p(x)$ is a radially symmetric (resp.\ axially symmetric with respect to some 
half-space) function, then some symmetry (resp.\ partial symmetry) results where recently obtained in \cite{MSS}
for semi-stable solutions and mountain-pass solutions.

\bigskip

\end{document}